\documentclass{amsart}

\newtheorem{theorem}{Theorem}[section]
\newtheorem{lemma}[theorem]{Lemma}

\theoremstyle{definition}
\newtheorem{definition}[theorem]{Definition}

\newtheorem{note}[theorem]{Note}
\newtheorem{proposition}[theorem]{Proposition}
\theoremstyle{remark}
\newtheorem{remark}[theorem]{Remark}
\newtheorem{corollary}[theorem]{Corollary}

\numberwithin{equation}{section}

\begin{document}

\title{Dimensional analysis of fractal interpolation functions}



\author{S. Verma}
\address{Department of Mathematics, IIT Delhi, New Delhi, India 110016 }
\email{saurabh331146@gmail.com}
\author{S. Jha}
\address{Department of Mathematics, NIT Rourkela, India 769008 }
\email{jhasa@nitrkl.ac.in}



 
\begin{abstract}
We provide a rigorous study on dimensions of fractal interpolation function defined on a closed and bounded interval of $\mathbb{R}$ which is associated to a continuous function with respect to a base function, scaling functions and a partition of the interval. In particular, we provide an exact estimation of the box dimension of  $\alpha$-fractal functions.

\end{abstract}

\maketitle


.

\section{INTRODUCTION}
The idea of fractal interpolation functions was introduced by Barnsley \cite{MF1}. Many authors attempted to calculate the box and Hausdorff dimensions of the graph of fractal interpolation functions corresponding to a data set, see for instance \cite{MF1,Hardin,M2,RSY,Vijender2}. Furthermore, some authors \cite{M1,WY} 
 have studied the properties of fractal interpolation functions corresponding to a data set.

Here we first start with iterated function system, for more details see \cite{MF2}.
\subsection{Iterated Function System} Let $(X, d)$ be a complete metric space, and let $H(X)$ be the family of all
nonempty compact subsets of $X$. We define the Hausdorff metric
$$h(A,B):=  \max \Big\{\max_{a \in A} \min_{b \in B} d(a,b),  \max_{b \in B} \min_{a \in A} d(a,b)  \Big\}.$$
It is well known that $(H(X),h)$ is a complete metric space. Let $k$ be a positive integer, and let, for $i=1, 2, \dots k$, $w_i$
be contraction selfmap of $X$, i.e., there exist  real numbers $R_i \in [0,1)$ such that
$$d(w_i(x),w_i(y)) \le  R_i d(x,y) \quad \forall \quad x , y \in X.$$
\begin{definition}
The system $\big\{(X,d); w_1,w_2,\dots,w_k \big\}$ is called an iterated function system, IFS for short.
\end{definition}
The IFS generates the mapping $W$ from $H(X)$ into $H(X)$ given by
$$ W(A) = \cup_{i=1}^k w_i (A).$$
The Hutchinson-Barnsley map $W$ defined above is then a contraction
mapping, with respect to the Hausdorff metric $h$, the Lipschitz
constant $R_*:=\max\{R_1, R_2, \dots, R_k\}$.
  Thus, by the Banach
contraction principle, there exists a unique nonempty compact subset $A$ such that
$A= \cup_{i=1}^k w_i(A)$. Such a set $A$ is
termed the attractor of the IFS.\\
The reader is referred to \cite{MF1,MF2,N2} for the upcoming technical introduction. The method of constructing fractal interpolation functions (FIFs) is as follows:
\subsection{Fractal Interpolation Functions}
 Consider a set of interpolation points $ \{(x_n,y_n) : n=1,2,\dots,N\} $ with strictly increasing abscissa. Set $J = \{1,2,...,N-1\}$,  $I= [x_1, x_N] $ and for $ j \in J,$  let $I_j= [x_j, x_{j+1}]$.  For $ j\in J$, let $L_j:I \rightarrow I_j $ be a contraction homeomorphism such that
$$ L_j(x_1)=x_j, L_j(x_N)=x_{j+1}, j \in J.$$
 For $j\in J $, let $F_j: I\times \mathbb{R} \rightarrow \mathbb{R} $ be a mapping satisfying
$$ |F_j(x,y) - F_j(x,y_*)| \leq \kappa_j |y- y_*| ,$$
$$ F_j(x_1,y_1)=y_j, F_j(x_N,y_N)=y_{j+1}, j \in J ,$$
where $(x,y), (x,y_*) \in K $ and $ 0 \leq \kappa_j < 1 $ for all $j \in J .$ We shall take
$$L_j(x)=a_j x+ b_j, \quad F_j(x,y) = \alpha_j y + q_j (x).$$
In the above expressions $a_j$ and $b_j$ are determined so that the conditions $ L_j(x_1)=x_j, L_j(x_N)=x_{j+1}$ are satisfied. The multipliers  $\alpha_j$, called scaling factors,  are such that $-1< \alpha_j <1$ and $q_j:I\rightarrow \mathbb{R}, j \in J $ are suitable continuous functions satisfying the ``join-up conditions'' imposed for the bivariate maps $F_j$. That is, $q_j(x_1)=y_j-\alpha_j y_1$ and $q_j(x_N)=y_{j+1}-\alpha_j y_N$ for all $j \in J $.
Now define functions $W_j:I\times \mathbb{R} \rightarrow  I\times \mathbb{R} $  for $j \in J$ by  $$W_j(x,y)=\big(L_j(x),F_j(x,y)\big). $$
 Theorem $1$ in \cite{MF1} says that the IFS $\mathcal{I}:=\{I\times \mathbb{R};W_1,W_2,\dots,W_{N-1}\}$ defined above has a unique attractor which is the graph of a function $g$ which satisfies
the following functional equation reflecting self-referentiality:
     $$g(x)= \alpha_j g \big(L_j^{-1}(x) \big)+ q_j \big(L_j^{-1}(x)\big), x \in I_j, j \in J.$$
   In \cite{MF1}, Barnsley gave an estimate for the Hausdorff dimension of an \emph{affine} FIF. Falconer \cite{Fal} also estimated the Hausdorff dimension of an \emph{affine} FIF. Further, in \cite{MF2,MF6,Hardin}, Barnsley and his collaborators calculated the box dimensions for classes of \emph{affine} FIFs. In \cite{MF3}, Barnsley and Massopust computed the box dimensions of FIFs generated by bilinear maps. In \cite{HM}, Hardin and Massopust produced a formula for the box dimension of vector-valued multivariate FIFs. In this article, we focus on dimensions of a special type of FIFs known as $\alpha$-fractal function.
 \subsection{$\alpha$-Fractal Functions: a Fractal Perturbation Process}
 The idea in the construction of a FIF can be adapted to obtain a class of fractal functions associated with a prescribed continuous function on a compact interval in $\mathbb{R}$. To this end, as is customary, let us denote by $\mathcal{C}(I)$, the space of all continuous real-valued functions defined on a compact interval $I=[x_1,x_N]$ in $\mathbb{R}$. We shall endow $\mathcal{C}(I)$ with the uniform norm. Let $f \in \mathcal{C}(I)$ be prescribed, referred to as the germ function. Let us consider the following elements to construct the IFS.
  \begin{enumerate}
 \item  A partition $\Delta:=\{x_1,x_2,\dots,x_N:x_1<x_2<\dots<x_N\}$ of $I=[x_1,x_N]$.

 \item For each $j \in J$, let $\alpha_j: I \rightarrow \mathbb{R}$ be continuous functions with $\| \alpha_j\|_{\infty}=\max\limits_{j}\{|\alpha_j(x)|:x\in I\} < 1$. These functions are referred to as scaling functions. Consider  $\alpha= \big(\alpha_1, \alpha_2, \dots, \alpha_{N-1}\big) \in \big(\mathcal{C}(I) \big)^{N-1}$, referred to as scaling vector.

 \item A continuous function $b:I \rightarrow \mathbb{R}$ such that  $b(x_1)=f(x_1), b(x_1)= f(x_N)$ and $b \ne f$, termed base function.
 \end{enumerate}
 Let us define maps
  \begin{equation} \label{BE1}
\begin{split}
L_j(x) =&~ a_j x + b_j, \\
F_j(x,y)=&~ \alpha_j(x) y + f \big( L_j(x)\big) -\alpha_j(x) b(x).
\end{split}
\end{equation}
By Theorem $1$ in \cite{MF1} one can see that the corresponding IFS $\mathcal{I}:=\{I\times \mathbb{R};W_1,W_2,\dots,W_{N-1}\}$, where
$$W_j(x,y) = \Big(L_j(x), F_j(x,y)\Big),$$
has a unique attractor, which is the graph of a continuous function, denoted by  $f_{\Delta,b}^{\alpha}: I \rightarrow \mathbb{R}$ such that $f_{\Delta,b}^{\alpha}(x_n)= f(x_n), n=1,2,\dots,N $.  Furthermore, $f_{\Delta,b}^{\alpha}$ satisfies the self-referential equation
$$ f_{\Delta,b}^{\alpha}(x)= f(x)+\alpha_j(L_j^{-1}(x)).(f^{\alpha}- b)\big(L_j^{-1}(x)\big)~~~~\forall~~ x \in I_j,~~ j \in J.$$
  It is observed that generally $f_{\Delta,b}^{\alpha}$ is non-differentiable (for less restrictions on $ \alpha$ and $b$) and its graph has a non-integer Hausdorff-Besicovitch dimension. One can treat  $f_{\Delta,b}^{\alpha}$ as a ``fractal perturbation" of $f$.
 \begin{definition}
 The function  $f_{\Delta,b}^{\alpha}$ is called the $\alpha$-fractal function corresponding to $f$ with respect to $\Delta$ and $b$.
 \end{definition}
There has been a great interest to study the properties of $\alpha$-fractal function $f^{\alpha}$, the reader is referred to \cite{N2,N1}.
 Akhtar et-al \cite{AGN2} calculated the box dimension of graph of $f^{\alpha}.$ Recently the authors of \cite{VV} estimated the box dimension of the graph of $f^{\alpha}$. We refer the reader to \cite{VV1} for the bivariate $\alpha$-fractal functions and further developments.
In this article, we study the box and Hausdorff dimensions of graph of $f^{\alpha}$ deeply. Our results are generalizing certain existing results in an exciting manner.
 
We skip definitions of the box and Hausdorff dimensions and refer the reader to \cite{Fal} for their definitions.

\section{Main Theorems}
We write the following note which is a modification of \cite[Proposition $1$]{MF3}.
\begin{note}
We define a metric $d$ on $I \times \mathbb{R}$ as follows $$d((x,y),(z,w))= c_1|x-z| +c_2|(y-f^{\alpha}(x))-(w-f^{\alpha}(z))|~ \forall (x,y),(z,w) \in I \times \mathbb{R}.$$ Then $d$ is a metric on $I \times \mathbb{R}.$ Furthermore, $(I \times \mathbb{R};d)$ is a complete metric space.
\end{note}
\begin{remark}
We recall the functional equation
$$ f^{\alpha}(x)= f(x)+\alpha_j(L_j^{-1}(x)).(f^{\alpha}- b)\big(L_j^{-1}(x)\big)~~~~\forall~~ x \in I_j,~~ j \in J.$$
Now for every $ x \in I_j$ and $ j \in J$, we have 
\begin{equation*}
      \begin{aligned}
           |f^{\alpha}(x)- f(x)| & =|\alpha_j(L_j^{-1}(x)).(f^{\alpha}- b)\big(L_j^{-1}(x)\big)|\\ & =|\alpha_j(L_j^{-1}(x))|~|(f^{\alpha}- b)\big(L_j^{-1}(x)\big)|\\ &  \le \|\alpha_j\|_{\infty} \|f^{\alpha}- b\|_{\infty} \\&  \le \|\alpha\|_{\infty} \|f^{\alpha}- b\|_{\infty}.
      \end{aligned}
     \end{equation*}  
     The above implies that $\|f^{\alpha}- f\|_{\infty} \le \|\alpha\|_{\infty} \|f^{\alpha}- b\|_{\infty}.$ Using triangle inequality we obtain $\|f^{\alpha}- f\|_{\infty} \le \|\alpha\|_{\infty} \|f^{\alpha}- f\|_{\infty}+\|\alpha\|_{\infty} \| f -b \|_{\infty}.$ Finally we have $\|f^{\alpha}- f\|_{\infty} \le \frac{\|\alpha\|_{\infty}}{1-\|\alpha\|_{\infty}}\| f -b \|_{\infty}.$
Further, we have $$\|f^{\alpha}\|_{\infty} - \| f\|_{\infty}\le \|f^{\alpha}- f\|_{\infty} \le \frac{\|\alpha\|_{\infty}}{1-\|\alpha\|_{\infty}}\| f -b \|_{\infty}.$$ Therefore, we get $ \|f^{\alpha}\|_{\infty} \le \|f\|_{\infty} +\frac{\|\alpha\|_{\infty}}{1 - \|\alpha \|_{\infty}} \|f-b\|_{\infty}:=M.$
\end{remark}
We should mention that the techniques involved in the proof of the next proposition is same as \cite[Theorem $4$]{MF3}.
\begin{proposition}\label{prop1}
The map $W_j: I \times [-M, M] \to I \times [-M, M] $ is a contraction map with respect to the above metric provided $$\max \Big\{a_j+ \frac{2c_2Mk_{\alpha_j}}{c_1}, \| \alpha_j \|_{\infty} \Big\} < 1$$
and $\alpha_j :I \to \mathbb{R} $ satisfies $|\alpha_j(x)- \alpha_j(y)| \le k_{\alpha_j} |x-y|.$
\end{proposition}
\begin{proof}
Let $(x,y),(z,w) \in I \times [-M, M].$ We have 
\begin{equation*}
\begin{aligned}
& d\big(W_j(x,y),W_j(z,w)\big)-c_1|L_j(x)-L_j(z)|\\= & ~ c_2\Big|\Big(\alpha_j(x) y + f \big( L_j(x)\big) -\alpha_j(x) b(x)-f^{\alpha}\big( L_j(x)\big)\Big)\\&-\Big(\alpha_j(z) w + f \big( L_j(z)\big) -\alpha_j(z) b(z)-f^{\alpha}\big( L_j(z)\big)\Big)\Big|\\ =& ~ c_2 \Big| \alpha_j(x)\big(y-f^{\alpha}(x)\big)- \alpha_j(z)\big(w-f^{\alpha}(z)\big)\Big|\\ = & ~ c_2 \Big| \alpha_j(x)\big(y-f^{\alpha}(x)\big)-\alpha_j(x)\big(w-f^{\alpha}(z)\big) + \alpha_j(x)\big(w-f^{\alpha}(z)\big)-\alpha_j(z)\big(w-f^{\alpha}(z)\big)\Big|\\ \le & ~c_2 |\alpha_j(x) |~ \big|\big(y-f^{\alpha}(x)\big)-\big(w-f^{\alpha}(z)\big)\big|\\& + c_2 |\alpha_j(x) - \alpha_j(z)|~ |w- f^{\alpha}(z) | \\ \le & ~c_2 \| \alpha_j \|_{\infty} | \big(y-f^{\alpha}(x)\big)-\big(w-f^{\alpha}(z)\big)|+ 2Mc_2 k_{\alpha_j} |x-z|.
\end{aligned}
\end{equation*}
Further, we get
\begin{equation*}
\begin{aligned}
d\big(W_j(x,y),W_j(z,w)\big)  \le & ~c_1 |L_j(x)-L_j(z)|\\ & + c_2 \| \alpha_j \|_{\infty} | \big(y-f^{\alpha}(x)\big)-\big(w-f^{\alpha}(z)\big)|+ 2Mc_2 k_{\alpha_j} |x-z|\\ = & ~c_1 a_j |x -z| \\ & + c_2 \| \alpha_j \|_{\infty} | \big(y-f^{\alpha}(x)\big)-\big(w-f^{\alpha}(z)\big)|+ 2M c_2k_{\alpha_j} |x-z|\\ = & \Big(a_j+ \frac{2c_2Mk_{\alpha_j}}{c_1}\Big)c_1|x-z| \\ &+ c_2 \| \alpha_j \|_{\infty} | \big(y-f^{\alpha}(x)\big)-\big(w-f^{\alpha}(z)\big)|\\ \le & \max \Big\{a_j+ \frac{2c_2Mk_{\alpha_j}}{c_1}, \| \alpha_j \|_{\infty} \Big\} \Big(c_1|x-z| \\& +c_2\big(y-f^{\alpha}(x)\big)-\big(w-f^{\alpha}(z)\big)| \Big)\\ =&
\max \Big\{a_j+ \frac{2c_2Mk_{\alpha_j}}{c_1}, \| \alpha_j \|_{\infty} \Big\}~ d\big((x,y),(z,w)\big). 
\end{aligned}
\end{equation*}
This completes the proof.
\end{proof}

We say that an IFS $\{X;w_1,w_2,\dots,w_k\}$ satisfies the open set condition if there exists a non-empty open set $U$ with $$ \bigcup\limits_{i=1}^k w_i(U) \subset U$$ and terms present in the above union are disjoint. Further, if the above $U$ satisfies $U \cap A \ne \emptyset $ then we call that the IFS satisfies the strong open set condition.
Now we are ready to prove the following.
 
\begin{theorem}
Let $\mathcal{I}:=\{I\times \mathbb{R};W_1,W_2,\dots,W_{N-1}\}$ be the IFS as defined earlier such that $$  r_i \|(x,y)-(w,z) \|_2 \le \|W_j(x,y) - W_j(w,z)\|_2 \le R_i \|(x,y)-(w,z) \|_2 ,$$ for every $(x,y),(w,z) \in I\times \mathbb{R},$ where $0 < r_i \le R_i < 1 ~ \forall~ i \in \{1,2,\dots ,N-1\}.$ Then $s_* \le  \dim_H(Graph(f^{\alpha})) \le s^*  ,$ where $s_*$ and $s^*$ are determined by  $ \sum\limits_{i=1}^{N} r_i^{s_*} =1$ and $ \sum\limits_{i=1}^{N} R_i^{s^*} =1$ respectively.
\end{theorem}

\begin{proof}
Following Proposition $9.6$ in \cite{Fal} we have the required upper bound. For the lower bound of Hausdorff dimension of $Graph(f^{\alpha})$ we proceed as follows.

Define $U= (x_1,x_N) \times \mathbb{R}.$ It is plain to see that $$L_i\big((x_1,x_N)\big) \cap L_j\big((x_1,x_N)\big)=\emptyset,$$ for every $i,j \in J:=\{1,2,\dots , N-1\}$ with $i \ne j.$ This immediately yields $$W_i(U) \cap W_{j}(U)=\emptyset,$$ for every $i,j \in J$ satisfying $ i \ne j.$ Using $U \cap Graph(f^{\alpha}) \ne \emptyset,$ one deduces that the IFS satisfies the SOSC.
Since  $U \cap Graph(f^{\alpha}) \ne \emptyset $, we have an index $i \in J^*$ with $Graph(f^{\alpha})_i \subset U$, where $J^*:=\cup_{m \in \mathbb{N}}\{1,2,\dots, N-1\}^m$, that is, the set of all finite sequences which are made up of the elements of $J$ and $Graph(f^{\alpha})_i:=W_i(Graph(f^{\alpha})):=W_{i_1} \circ W_{i_2}\circ \dots \circ W_{i_m} (Graph(f^{\alpha}))$ for $i \in J^m$($m$-times Cartesian product of $J$ with itself) and for $m \in \mathbb{N}.$ Now, it is obvious that for any $n$ and $j \in J^n$, the sets $Graph(f^{\alpha})_{ji}$ are disjoint. Furthermore, the IFS $\{W_{ji}: j \in J^n\}$ satisfies the hypotheses of Proposition $9.7$ in \cite{Fal}. Therefore, with the notation $r_j=r_{j_1}r_{j_2} \dots r_{j_n}$ we have $ s_n \le \dim_H(G^*)$, where $G^* $ is an attractor of the aforesaid IFS and $ \sum_{ j \in J^n} r_{ji}^{s_n} =1.$ Since $G^* \subset Graph(f^{\alpha})$, we get $  s_n \le \dim_H(G^*) \le \dim_H(Graph(f^{\alpha}))$. Suppose that $ \dim_H(Graph(f^{\alpha})) < s_*,$ where $ \sum_{i=1}^{N-1} r_{i}^{s_*} =1.$ This gives $ s_n < s_*.$ Using the previous estimates, we have
\begin{equation}
\begin{aligned}
r_{i}^{- s_n} &= \sum_{ j \in J^n} r_{j}^{s_n}\\ & \ge \sum_{ j \in J^n} r_{j}^{\dim_H(Graph(f^{\alpha}))}\\& = \sum_{ j \in J^n} r_{j}^{s_*} r_{j}^{\dim_H(Graph(f^{\alpha})) -s_*}\\& \ge \sum_{ j \in J^n} r_{j}^{s_*} r_{max}^{n(\dim_H(Graph(f^{\alpha})) - s_*)}\\& = r_{max}^{n(\dim_H(Graph(f^{\alpha})) -s_*)},
\end{aligned}
\end{equation}
where $r_{max}=\max\{r_1, r_2, \dots,r_{N-1}\}.$ 
Since $ r_{max} < 1$ and the term on left side in the above expression is bounded, we have a contradiction as $n $ tends to infinity. Thus our claim is wrong. This implies that $ \dim_H(Graph(f^{\alpha})) \ge  s_*,$ which is the required result.
\end{proof}

\begin{remark}
Under the assumptions of Proposition \ref{prop1} we may find a upper bound for the Hausdorff dimension of graph of $\alpha$-fractal function $f^{\alpha}$ using the above theorem.
\end{remark}

\begin{remark}
Note that the open set $(0,1) \times \mathbb{R}$ will serve for our purpose to satisfy the strong open set condition for $\{I \times \mathbb{R};W_j: j=1,2 \dots, N-1\}.$ With the aid of the above theorem we are able to improve Theorem $4$ in \cite{MF1}. In particular, with the notation in \cite{MF1} we can omit the following condition from that theorem $$ t_1.t_N \le (Min\{a_1,a_N\})\Big(\sum_{n=1}^{N}t_n^l\Big)^{2/l}.$$
\end{remark}

\begin{remark}
Schief \cite{Schief} followed the technique of Bandt \cite{Bandt} and proved that open set and strong open set conditions are equivalent for similitudes. Further Peres, Rams, Simon, and Solomyak \cite{Peres} showed Schief's theorem for self-conformal maps. The same result with a different approach can be seen in \cite{Lau1,Lau2, Lau3}. We do not know whether or not the open set and strong open set conditions are equivalent for $\alpha$-fractal functions.
\end{remark}
\begin{remark}
It is known that for a pure self-similar set or self-conformal set $A,\;  \dim_H(A)=\underline{\dim}_B(A)= \overline{\dim}_B(A),$ for more details see \cite{Fal2}. Note that  the nature of $\alpha$-fractal functions depend on the IFS parameter. In particular, one can obtain pure self-similar or partial self similar $\alpha$-fractal functions by choosing suitable scaling functions and thus for $\alpha$-fractal functions we may or may not get the equal dimensions. 
\end{remark}
\begin{remark}
Here we talk about continuity of the Hausdorff dimension. Note that, in general Hausdorff dimension is not a continuous function. For example, $A_n:=[-\frac{1}{n},\frac{1}{n}] \to A:=\{0\}$ in Hausdorff metric but $\dim_H(A_n)=1 $ does not converge to $\dim_H(A)=0.$ In \cite{VV}, the continuous dependence of $\alpha$-fractal function on the parameters is studied. One may pose a question of continuity of the Hausdorff dimension of $\alpha$-fractal function with respect to the parameters involved. However, it seems that the result may not hold in general.
\end{remark}

\section{Oscillation Spaces}
We refer the reader to \cite{ACar,ADBJ} for oscillation spaces.
Let $Q \subset [0,1]$ $p$-adic subinterval so that $ Q =\Big[\frac{i}{p^m},\frac{i+1}{p^m}\Big]$ for some integers $m \ge 0$ and $0 \le i < \frac{1}{p^m}.$
For a continuous function $f:[0,1] \to \mathbb{R}$, we define oscillation of $f$ over $Q$ 
\begin{equation*}
\begin{aligned}
R_f(Q)=& \sup_{x,y \in Q} |f(x) -f(y)|\\=& \sup_{x \in Q}f(x) - \inf_{y \in Q} f(y),
\end{aligned}
\end{equation*}
and total oscillation of order $m$, $$ Osc(m,f)= \sum_{|Q|=p^{-m}} R_f(Q),$$ where the sum ranges over all $p$-adic intervals $ Q \subset [0,1]$ of length $|Q|=\frac{1}{p^m}.$\\
 Let $\beta\in \mathbb{R}$. We define the oscillation space $\mathcal{V}^{\beta}(I)$ by $$ \mathcal{V}^{\beta}(I)= \Big\{f \in \mathcal{C}(I): \sup_{m \in \mathbb{N} }\frac{ Osc(m,f)}{p^{m(1-\beta)}}< \infty \Big\}.$$
One can define $$ \mathcal{V}^{\beta-}(I)= \{f \in \mathcal{C}(I): f \in \mathcal{V}^{\beta- \epsilon}(I)~ \forall ~ \epsilon >0 \},$$ and 
$$ \mathcal{V}^{\beta+}(I)= \{f \in \mathcal{C}(I): f \notin \mathcal{V}^{\beta+ \epsilon}(I)~ \forall ~ \epsilon >0 \}.$$
Now, let us write the next two theorems as a prelude.
\begin{theorem}[\cite{ACar}, Theorem $4.1$]\label{useit77}
Let $f$ be a real-valued continuous function defined on $I$, we have
$$\overline{\dim}_B(Graph(f)) \le 2 - \gamma \iff f \in \mathcal{V}^{\gamma-}(I) ~~ ~\text{if}~ 0 <\gamma \le 1$$
and 
$$\overline{\dim}_B(Graph(f)) \ge 2 - \gamma \iff f \in \mathcal{V}^{\gamma+}(I) ~ ~\text{if}~ 0 \le \gamma < 1.$$
\end{theorem}
\begin{theorem}[\cite{ADBJ}, Theorem $3.1$]\label{useit78}
Let $f:I \to \mathbb{R}$ be a continuous function and let $0 < \gamma < 1.$ Then we have
\[\overline{\dim}_B(Graph(f)) =  2 - \gamma \iff f \in \cap_{\theta < \gamma} \mathcal{V}^{\theta}(I) \backslash \cup_{\beta > \gamma} \mathcal{V}^{\beta}(I).\]
\end{theorem}

\begin{lemma}\label{Lem1}
Let $f,g \in \mathcal{C}(I)$ and $\lambda \in \mathbb{R}.$ Then, for $m \in \mathbb{N}$ we have the following
\begin{enumerate} 
\item $Osc(m,\lambda f)=| \lambda | Osc(m,f)$
\item  $Osc(m, f+g) \le  Osc(m,f)+ Osc(m,g)$
\item  $Osc(m,fg) \le  \|g\|_{\infty} Osc(m,f)+\|f\|_{\infty} Osc(m,g).$
\end{enumerate}
\end{lemma}

\begin{proof}
\begin{enumerate}
\item 
For $m \in \mathbb{N}$ and $f,g \in \mathcal{C}(I)$, one proceeds as follows
\begin{equation*}
\begin{aligned}
Osc(m,\lambda f) = & \sum_{|Q|=p^{-m}} R_{\lambda f} (Q)\\ =&  \sum_{|Q|=p^{-m}} \sup_{x,y \in Q} |(\lambda f)(x)-(\lambda f)(y)|\\ =&
\sum_{|Q|=p^{-m}}|\lambda | \sup_{x,y \in Q} |f(x)- f(y)|\\ =&
|\lambda |\sum_{|Q|=p^{-m}} \sup_{x,y \in Q} | f(x)-f(y)|\\ =& | \lambda | ~ Osc(m,f).
\end{aligned}
\end{equation*}
\item Turning to second item we have
\begin{equation*}
\begin{aligned}
Osc(m, f+g) = & \sum_{|Q|=p^{-m}} R_{ f+g} (Q)\\ =&  \sum_{|Q|=p^{-m}} \sup_{x,y \in Q} |( f+g)(x)-( f +g )(y)|\\ \le &
\sum_{|Q|=p^{-m}}\Big( \sup_{x,y \in Q} |f(x)- f(y)|+\sup_{x,y \in Q} |g(x)- g(y)|\Big)\\ =&
\sum_{|Q|=p^{-m}}R_{ f} (Q)+\sum_{|Q|=p^{-m}}R_{ g} (Q)\\ = & ~Osc(m,f) +Osc(m,g).
\end{aligned}
\end{equation*}
\item The third item follows through the following lines.
\begin{equation*}
\begin{aligned}
Osc(m, fg) = & \sum_{|Q|=p^{-m}} R_{ fg} (Q)\\ =&  \sum_{|Q|=p^{-m}} \sup_{x,y \in Q} |( fg)(x)-( f g )(y)|\\ = &  \sum_{|Q|=p^{-m}} \sup_{x,y \in Q} |f(x)g(x)-f(y)g(x) +f(y)g(x)-f (y)g (y)| \\ \le &
\sum_{|Q|=p^{-m}}\Big( \sup_{x,y \in Q} |g(x)||f(x)- f(y)|+\sup_{x,y \in Q} |f(y)||g(x)- g(y)|\Big)\\ \le &
\sum_{|Q|=p^{-m}}\Big(\|g\|_{\infty} \sup_{x,y \in Q}|f(x)- f(y)|+\|f\|_{\infty}\sup_{x,y \in Q} |g(x)- g(y)|\Big)\\ =&
\|g\|_{\infty}\sum_{|Q|=p^{-m}}R_{ f} (Q)+\|f\|_{\infty}\sum_{|Q|=p^{-m}}R_{ g} (Q)\\ =& \|g\|_{\infty}Osc(m,f) +\|f\|_{\infty}Osc(m,g).
\end{aligned}
\end{equation*}
\end{enumerate}
Thus, the proof of the lemma is complete.
\end{proof}

\begin{proposition}\label{prop2}
Let $f \in \mathcal{V}^{\beta}(I),$ we define $\|f\|_{\mathcal{V}^{\beta}}:= \|f\|_{\infty} + \sup_{m \in \mathbb{N} }\frac{ Osc(m,f)}{p^{m(1-\beta)}}.$ Then $\|.\|_{\mathcal{V}^{\beta}}$ forms a norm on $\mathcal{V}^{\beta}(I).$
\end{proposition}
\begin{proof}
\begin{enumerate}
Through simple and straightforward calculations, we have

\item 
\begin{equation*}
\begin{aligned}
& \|f\|_{\mathcal{V}^{\beta}}= 0 \\ \iff & \|f\|_{\infty}=0 ~\text{and}~~ \sup_{m \in \mathbb{N} }\frac{ Osc(m,f)}{p^{m(1-\beta)}}=0 \\ \iff & f = 0,
\end{aligned}
\end{equation*}

\item
\begin{equation*}
\begin{aligned}
 \|\lambda f\|_{\mathcal{V}^{\beta}}= & \|\lambda f\|_{\infty} + \sup_{m \in \mathbb{N} }\frac{ Osc(m,\lambda f)}{p^{m(1-\beta)}} \\ = &  |\lambda|~ \|f\|_{\infty} + |\lambda| \sup_{m \in \mathbb{N} }\frac{ Osc(m,f)}{p^{m(1-\beta)}} \\ = & | \lambda |~ \|f\|_{\mathcal{V}^{\beta}},
\end{aligned}
\end{equation*}
and
\item
\begin{equation*}
\begin{aligned}
 \| f + g\|_{\mathcal{V}^{\beta}}= & \| f + g\|_{\infty} + \sup_{m \in \mathbb{N} }\frac{ Osc(m, f+g)}{p^{m(1-\beta)}} \\ \le  &   \|f\|_{\infty} +\|g \|_{\infty} +  \sup_{m \in \mathbb{N} }\frac{ Osc(m,f+g)}{p^{m(1-\beta)}} \\ \le  &   \|f\|_{\infty} +\|g \|_{\infty} +  \sup_{m \in \mathbb{N} }\frac{ Osc(m,f)}{p^{m(1-\beta)}} +\sup_{m \in \mathbb{N} }\frac{ Osc(m,g)}{p^{m(1-\beta)}} \\ = &  \|f\|_{\mathcal{V}^{\beta}}+\|g\|_{\mathcal{V}^{\beta}},
\end{aligned}
\end{equation*}
\end{enumerate}
hence the proof.
\end{proof}
\begin{lemma}\label{lem2}
Let $(f_n)$ be a sequence of continuous functions that converges uniformly to some $f: I \to \mathbb{R}$ and $m \in \mathbb{N},$ then we have $$ Osc(m,f_n) \to Osc(m,f).$$
Furthermore, let $(f_n)$ be a sequence in $\mathcal{V}^{\beta}(I)$ that converges uniformly to some $ f:I \to \mathbb{R},$ then we have $$ \sup_{m \in \mathbb{N}} \frac{Osc(m,f)}{p^{m(1-\beta)}} \le \liminf_{n \to \infty} \sup_{m \in \mathbb{N}} \frac{Osc(m,f_n)}{p^{m(1-\beta)}}.$$
\end{lemma}
\begin{proof}
Let $m \in \mathbb{N}$, we have
\begin{equation*}
\begin{aligned}
\lim_{n \to \infty} Osc(m,f_n)= & \lim_{n \to \infty} \sum_{|Q|=p^{-m}} R_{f_n} (Q)\\ =& \lim_{n \to \infty} \sum_{|Q|=p^{-m}} \sup\limits_{x,y \in Q} |f_n(x)-f_n(y)|\\ =& \sum_{|Q|=p^{-m}} \lim_{n \to \infty}  \sup\limits_{x,y \in Q} |f_n(x)-f_n(y)|\\ =&  \sum_{|Q|=p^{-m}} \sup\limits_{x,y \in Q} |f(x)-f(y)|\\ =&~ Osc(m,f).
\end{aligned}
\end{equation*}
Now for $m \in \mathbb{N},$ we get 
\begin{equation*}
\begin{aligned}
\frac{Osc(m,f)}{p^{m(1-\beta)}} = & \frac{\sum\limits_{|Q|=p^{-m}} \sup\limits_{x,y \in Q} |f(x)-f(y)|}{p^{m(1-\beta)}}\\ = & \frac{\sum\limits_{|Q|=p^{-m}} \sup\limits_{x,y \in Q} |\lim\limits_{n \to \infty} f_n(x)-\lim\limits_{n \to \infty}f_n(y)|}{p^{m(1-\beta)}}\\ = & \lim_{n \to \infty} \frac{\sum\limits_{|Q|=p^{-m}} \sup\limits_{x,y \in Q} | f_n(x)-f_n(y)|}{p^{m(1-\beta)}}\\ = &   \lim_{n \to \infty} \frac{Osc(m,f_n)}{p^{m(1-\beta)}}\\ \le  & 
\liminf_{n \to \infty} \left(\sup_{m \in \mathbb{N}}\frac{Osc(m,f_n)}{p^{m(1-\beta)}}\right).
\end{aligned}
\end{equation*}
Thus, the proof of the lemma is complete.
\end{proof}
\begin{theorem}
The space $(\mathcal{V}^{\beta}(I), \|.\|_{\mathcal{V}^{\beta}})$ is a Banach space.
\end{theorem}
\begin{proof}
Let $(f_n)$ be a Cauchy sequence in $\mathcal{V}^{\beta}(I)$ with respect to $\|.\|_{\mathcal{V}^{\beta}}.$ Equivalently, for a given $\epsilon >0$, we have a natural number $n_0$ such that $$ \|f_n -f_k \|_{\mathcal{V}^{\beta}}  < \epsilon ~ ~ \forall ~ n,k \ge n_0.$$
By definition of $\|.\|_{\mathcal{V}^{\beta}}$ one gets $ \|f_n -f_k \|_{\infty}  < \epsilon ~ ~ \forall` n,k \ge n_0.$ Since $(\mathcal{C}(I), \|.\|_{\infty})$ is a Banach space, we have a continuous function $f$ with $ \|f_n -f\|_{\infty} \to 0$ as $n \to \infty.$
We claim that $f \in \mathcal{V}^{\beta}(I)$ and $\|f_n -f\|_{V^{\beta}} \to 0$ as $n \to \infty.$ Let $m \in \mathbb{N}$ and $n \ge n_0.$ In view of Lemma \ref{lem2} we have 
\begin{equation*}
\begin{aligned}
\|f_n -f\|_{\infty} + \frac{Osc(m,f_n-f)}{p^{m(1-\beta)}} = & \lim_{k \to \infty} \Bigg(\|f_n -f_k\|_{\infty} + \frac{Osc(m,f_n-f_k)}{p^{m(1-\beta)}}\Bigg)\\  \le &  \lim_{k \to \infty} \Bigg(\|f_n -f_k\|_{\infty} + \sup_{m' \in \mathbb{N}}\frac{Osc(m',f_n-f_k)}{p^{m'(1-\beta)}}\Bigg)\\  \le &  \sup_{k \ge n_0} \Bigg(\|f_n -f_k\|_{\infty} + \sup_{m' \in \mathbb{N}}\frac{Osc(m',f_n-f_k)}{p^{m'(1-\beta)}}\Bigg)\\ = &  \sup_{k \ge n_0} \|f_n -f_k\|_{\mathcal{V}^{\beta}} \\ \le &~ \epsilon.
 \end{aligned}
 \end{equation*}
The above is true for every $m \in \mathbb{N}.$ Therefore, we obtain $f-f_{n_0} \in \mathcal{V}^{\beta}(I)$. Using Lemma \ref{Lem1} we have $f=f-f_{n_0}+f_{n_0} \in \mathcal{V}^{\beta}(I)$, and $\|f_n -f\|_{\mathcal{V}^{\beta}} \le \epsilon ~ \forall ~n \ge n_0$  
\end{proof}
\begin{remark}\label{reem}
If $|L_j(I)|= \frac{1}{p^{k_j}}$ for some $k_j \in  \mathbb{N}$ with $ \sum\limits_{j \in J} \frac{1}{p^{k_j}} = 1 $, then for $m \ge \max\limits_{j \in J} \{k_j\}$ we have $$Osc(m,f) = \sum\limits_{j \in J} Osc(m,f,L_j(I)),$$ where $Osc(m,f,L_j(I))= \sum\limits_{|Q|=p^{-m}, Q \subseteq L_j(I)} R_f(Q).$
\end{remark}
\begin{proof}
Since $I = \cup_{j \in J} L_j(I) $ and $ \sum\limits_{j \in J} \frac{1}{p^{k_j}} = 1 $ we have 
\begin{equation}
\begin{aligned}
Osc(m,f) & = \sum_{|Q|=p^{-m}} R_f(Q) \\
& = \sum_{j \in J} ~ ~\sum_{|Q|=p^{-m}, Q \subseteq L_j(I)} R_f(Q) \\ & 
= \sum_{j \in J} Osc(m,f,L_j(I)).
\end{aligned}
\end{equation}
\end{proof}
\begin{theorem}\label{BBVL3}
    Let $f ,b, \alpha_j~(j \in J) \in \mathcal{V}^{\beta}(I)$ be such that $b(x_1)=f(x_1)$ and $b(x_N)=f(x_N).$ Further we assume that $|L_j(I)|= \frac{1}{p^{k_j}}$ for some $k_j \in  \mathbb{N}$ with $ \sum\limits_{j \in J} \frac{1}{p^{k_j}} = 1 .$ For $\max \Bigg\{\| \alpha \|_{\infty} +\sum\limits_{j \in J} \sup_{m \in \mathbb{N}}\frac{Osc(m,\alpha_j)}{p^{m(1-\beta)}},\sum\limits_{j \in J}\|\alpha_j\|_{\infty}  \Bigg\}<1 $, we have $f^{\alpha} \in  \mathcal{V}^{\beta}(I) .$ 
 \end{theorem}
 \begin{proof}
     Let $ \mathcal{V}^{\beta}_f(I )= \{ g \in \mathcal{V}^{\beta}(I ): g(x_1)=f(x_1), ~g(x_N)=f(x_N) \}.$ We observe that the space $\mathcal{V}^{\beta}_f(I )$ is a closed subset of $\mathcal{V}^{\beta}(I).$ It follows that $\mathcal{V}^{\beta}_f(I )$ is a complete metric space with respect to the metric induced by norm $\|.\|_{\mathcal{V}^{\beta}}.$ We define a map $T: \mathcal{V}^{\beta}_f(I ) \rightarrow \mathcal{V}^{\beta}_f(I )$ by $$ (Tg)(x)=f(x)+\alpha_j(L_j^{-1}(x)) ~(g-b)(L_j^{-1}(x))  $$
     for all $x \in I_j $, where $j \in J.$
     First we observe that the mapping $T$ is well-defined.
     Using Remark \ref{reem}, for $g, h \in \mathcal{V}^{\beta}_f(I )$ we have 
          \begin{equation*}
               \begin{aligned}
                           \|Tg -Th\|_{\mathcal{V}^{\beta}} = & \|Tg -Th\|_{\infty} + \sup_{m \in \mathbb{N}}\frac{Osc(m,Tg-Th)}{p^{m(1-\beta)}}\\
                           \le &\| \alpha \|_{\infty} \|g -h\|_{\infty}+ \sum_{j \in J}\|\alpha_j\|_{\infty} \sup_{m \in \mathbb{N}}\frac{Osc(m,g-h)}{p^{m(1-\beta)}}\\ & + \sum_{j \in J}\|g -h\|_{\infty} \sup_{m \in \mathbb{N}}\frac{Osc(m,\alpha_j)}{p^{m(1-\beta)}}\\
                  \le & \Bigg(\| \alpha \|_{\infty} +\sum_{n \in J} \sup_{m \in \mathbb{N}}\frac{Osc(m,\alpha_j)}{p^{m(1-\beta)}} \Bigg) \| g-h\|_{\infty}\\ & + \Big(\sum_{j \in J}\|\alpha_j\|_{\infty}\Big) \sup_{m \in \mathbb{N}}\frac{Osc(m,g-h)}{p^{m(1-\beta)}}\\ \le & 
                  \max \Bigg\{\| \alpha \|_{\infty} +\sum_{j \in J} \sup_{m \in \mathbb{N}}\frac{Osc(m,\alpha_j)}{p^{m(1-\beta)}},\sum_{j \in J}\|\alpha_j\|_{\infty}  \Bigg\} \| g-h\|_{\mathcal{V}^{\beta}}.
               \end{aligned}
               \end{equation*}
       From the hypothesis, it follows that $T$ is a contraction map on $ \mathcal{V}^{\beta}_f(I ).$ Using the Banach contraction principle, we get a unique fixed point of $T$, namely $f^{\alpha} \in \mathcal{V}^{\beta}_f(I)$.
       Furthermore, since $T(f^{\alpha})=f^{\alpha}$, we write $f^{\alpha}$ as a part of the functional equation: $$ f^{\alpha}(L_j(x))= f(L_j(x))+\alpha_j(x) ~(f^{\alpha}-b)(x)$$ for every $x \in I$ and $j \in J.$ Now with $J:= \{1,2,3, \dots, N-1\}$ we define functions $W_j: I \times \mathbb{R} \rightarrow  I \times \mathbb{R} $  for $j \in J$ by $$W_j(x,y) = \Big(L_j(x), \alpha_j(x) y + f \big( L_j(x)\big) -\alpha_j(x) b(x)\Big).$$
       We show in the last part of the proof that graph of the associated fractal function $f^{\alpha}$ is an attractor of the IFS $\{I \times \mathbb{R}; W_j,j \in J\}.$ Following the proof of Theorem $1$ appeared in \cite{MF1} we may prove that attractor of the above IFS is graph of a function. It remains to show that it is actually the graph of fractal perturbation $f^{\alpha}$. To see that we use the functional equation and $I = \cup_{j \in J} L_j(I)$ and get
               \begin{equation*}
                          \begin{aligned} 
                             \cup_{j \in J} W_j(Graph(f^{\alpha}))= & \cup_{j \in J} \big\{W_j(x,f^{\alpha} (x)):x \in I \big\}\\
                             = & \cup_{j \in J} \Big\{\Big(L_j(x),\alpha_j(x) f^{\alpha} (x) + f \big( L_j(x)\big) -\alpha_j(x) b(x)\Big):x \in I \Big\}\\
                             = & \cup_{j \in J} \Big\{\big(L_j(x),f^{\alpha} (L_j(x))\big):x \in I \Big\}\\
                             = & \cup_{j \in J} \big\{(x,f^{\alpha} (x)):x \in L_j(I) \big\}\\
                             =& Graph(f^{\alpha}),
                             \end{aligned} 
                          \end{equation*}
completing the proof.
\end{proof}
\begin{remark}
Let $ 0<\gamma\leq 1$ and $f,b,\alpha_j$ be suitable functions satisfying the hypothesis of Theorem \ref{BBVL3}. Then, Theorem \ref{useit77} yields that $ \overline{\dim}_B(Graph(f^\alpha))\leq 2-\gamma$.
\end{remark}
\begin{remark}
Having Theorem \ref{useit78} in mind, we may ask the assumptions on the parameters for which $f^{\alpha} \in \cap_{\theta < \gamma} \mathcal{V}^{\theta}(I) \backslash \cup_{\beta > \gamma} \mathcal{V}^{\beta}(I).$ This question remains open.
\end{remark}
Before stating the upcoming remark, we define the H\"{o}lder space as follows:$$ \mathcal{H}^{s}(I ) := \{g:I \rightarrow \mathbb{R}: ~\text{g is H\"{o}lder continuous with exponent}~ s \} .$$
       If we equip the space $\mathcal{H}^{s}(I)$ with norm $ \|g\|_{\mathcal{H}}:= \|g\|_{\infty} +[g]_{s},$ where $$[g]_{s} = \sup_{x\ne y} \frac{|g(x)-g(y)|}{|x-y|^{s}}$$ then it forms a Banach space.
\begin{remark}
Let us start with the following example: the function $f:[0,1] \to \mathbb{R}$ defined by
$f(x) = |x - \frac{1}{2}|^{\beta},$ where $ 0 <\beta< 1,$ is a simple example of a function in $\cap_{\theta < 1} \mathcal{V}^{\theta}(I)$
which is only in the H\"{o}lder space $ \mathcal{H}^{\beta}(I ).$ Therefore, the dimension of the graph is $1$ while the classical
result only provides us with the upper bound $2 - \beta.$
Note that \cite{ADBJ} the spaces $\mathcal{V}^{\beta}(I)$ are refined version of H\"{o}lder spaces. Hence, our result obtained here generalizes many previous results, see, for instance, \cite{AGN2,VV}.
\end{remark}

The next theorem has been proved in \cite{VV} using the series expansion. We here give a different proof which we feel, is more general and direct.
\begin{theorem}\label{BBVL3}
    Let $f, b $ and $\alpha$ be  H\"{o}lder continuous with exponent $s$ such that $b(x_1)=f(x_1)$ and $b(x_N)=f(x_N).$ Then with the notation $a:= \min\{a_j: j \in J \}$ we have
     $f^{\alpha}$ is H\"{o}lder continuous with exponent $s$ provided $  \frac{\|\alpha\|_{\mathcal{H}}}{a^s}< 1 .$
 \end{theorem}
 \begin{proof}
     
     Let $ \mathcal{H}^{s}_f(I ):= \{ g \in \mathcal{H}^{s}(I ): g(x_1)=f(x_1), ~g(x_N)=f(x_N) \}.$ Applying the definition of closed set, we see that the set $\mathcal{H}^{s}_f(I )$ is a closed subset of $\mathcal{H}^{s}(I).$ Because $\mathcal{H}^{s}(I)$ is a Banach space as mentioned, it follows that $\mathcal{H}^{s}_f(I )$ is a complete metric space with respect to the metric induced by aforementioned norm $\|.\|_{\mathcal{H}}$ for $\mathcal{H}^{s}(I).$ We define a map $T: \mathcal{H}^{s}_f(I ) \rightarrow \mathcal{H}^{s}_f(I )$ by $$ (Tg)(x)=f(x)+\alpha_j(L_j^{-1}(x)) ~(g-b)(L_j^{-1}(x))  $$
     for all $x \in I_j $ where $j \in J.$
     First we shall show that $T$ is well-defined. For this let us note that
     \begin{equation*}
     \begin{split}
              [Tg]_s   = &\max_{j  \in J} \sup_{x \ne y, x,y  \in I_j} \frac{|Tg(x)-Tg(y)|}{|x-y|^{s}}\\
                 \le&  \max_{j  \in J} \Bigg[ \sup_{x \ne y, x,y \in I_j} \frac{|f(x)-f(y)|}{|x-y|^{s}}\\
                 & +  \sup_{x\ne y, x,y \in I_j} \frac{|\alpha_j(L_j^{-1}(x))| \Big|(g-b)(L_j^{-1}(x))-(g-b)(L_j^{-1}(y))\Big|}{|x-y|^{s}}\\ & + \sup_{x\ne y, x,y \in I_j} \frac{|(g-b)(L_j^{-1}(y))| \Big|\alpha_j(L_j^{-1}(x))-\alpha_j(L_j^{-1}(y))\Big|}{|x-y|^{s}}\Bigg]\\
                 \le & ~[f]_{s}+ \frac{\|\alpha\|_{\infty}}{a^{s}} \big( [g]_{s}+[b]_{s} \big)+ \frac{\|g-b\|_{\infty}}{a^{s}} [\alpha]_{s},
     \end{split}
     \end{equation*}
     where $[\alpha]_{s}= \max\limits_{j \in J} \sup\limits_{x \ne y, x,y \in I} \frac{ |\alpha_j(x)-\alpha_j(y)|}{|x-y|^{s}}.$
     For $g, h \in \mathcal{H}^{s}_f(I )$, we have
     \begin{equation*}
          \begin{aligned}
                      \|Tg -Th\|_{\mathcal{H}} &= \|Tg -Th\|_{\infty} + [Tg-Th]_{s}\\
                      &\le \| \alpha \|_{\infty} \|g -h\|_{\infty} + \frac{\|\alpha\|_{\infty}}{a^{s}} [g-h]_{s}+\frac{\|g-h\|_{\infty}}{a^{s}} [\alpha]_{s}\\
             &\le \frac{\|\alpha\|_{\mathcal{H}}}{a^{s}} \| g-h\|_{\mathcal{H}}.
          \end{aligned}
          \end{equation*}
  Since $\frac{\|\alpha\|_{\mathcal{H}}}{a^{s}} < 1$, it follows that $T$ is a contraction self map on $ \mathcal{H}^{s}_f(I ).$ Thanks to Banach contraction principle, a unique fixed point of $T$ exists. This proves the result.
    \end{proof}

\begin{remark}
 By dint of a more spirited effort (see \cite{VV}), we can observe with the help of the above theorem that the H\"{o}lder constant of the map $f^\alpha$ depends only on the germ function $f$, the partition $\Delta$ and the parameter maps $~b,~ \alpha$.
\end{remark}
Note that for equidistant nodes we have $a=a_j = \frac{1}{N-1}.$
\begin{theorem}\label{mainthm}
    Let $f$ be a germ function, and $b, \alpha_j$ be  suitable continuous functions such that
    \begin{equation}\label{Hypo}
    \begin{aligned}
    & |f(x) -f(y)| \le k_f |x-y |^{s},\\&
     |b(x) -b(y)| \le k_b |x-y|^{s},\\&
     |\alpha_j(x) -\alpha_j(y)| \le k_{\alpha} |x-y|^{s}
     \end{aligned}
     \end{equation}
      for every $x,y \in I ,j \in J,$ and for some $k_f, k_b, k_{\alpha} > 0, s\in (0,1]$.
      Further, assume that there are constants $K_f, \delta_0> 0$ such that for each $x \in I$ and $\delta < \delta_0$ there exists $y \in I$ with $|x-y| \le \delta$ and $ |f(x)-f(y)| \ge K_f |x -y|^{s}.$
     We have
           $\dim_B\big(Graph(f^{\alpha})\big) = 2 - s$ provided that $ \|\alpha\|_{\mathcal{H}}< a^s ~\min\Big\{1,\frac{K_f-(\|b\|_{\infty}+M)k_{\alpha}a^{-s} }{(k_{f^\alpha}+k_b)}\Big\}.$
\end{theorem}
\begin{proof}
In the light of Theorem \ref{BBVL3} and $ \|\alpha\|_{\mathcal{H}}< a^s$, we have $f^{\alpha}$ is H\"older continuous with the same exponent $s$. That is, we may consider  $$ |f^{\alpha}(x)-f^{\alpha}(y)| \le k_{f^{\alpha}} |x-y|^{s}$$ for some $ k_{f^{\alpha}}> 0. $ We obtain a bound for upper box dimension of the graph of the fractal function $f^{\alpha} $ as follows: For $0 < \delta < 1,$  let $N_{\delta}(Graph(f^{\alpha}))$ be the number of $\delta-$boxes that cover graph of $f^{\alpha},$ with $\left \lceil {.}\right \rceil $ the ceiling function, we have
\begin{equation}
\begin{aligned}
N_{\delta}(Graph(f^{\alpha})) & \le \sum_{i=1}^{\left \lceil {\frac{1}{\delta}}\right \rceil}  \Bigg(1+ \left \lceil { \frac{R_{f^{\alpha}}[(i-1)\delta,i\delta]}{\delta} }\right \rceil \Bigg) \\ & \le \sum_{i=1}^{\left \lceil {\frac{1}{\delta}}\right \rceil} \Bigg(2+  \frac{R_{f^{\alpha}}[(i-1)\delta,i\delta]}{\delta} \Bigg)\\& =  2 \left \lceil {\frac{1}{\delta}}\right \rceil + \sum_{i=1}^{\left \lceil {\frac{1}{\delta}}\right \rceil}   \frac{R_{f^{\alpha}}[(i-1)\delta,i\delta]}{\delta}\\& \le2 \left \lceil {\frac{1}{\delta}}\right \rceil + \sum_{i=1}^{\left \lceil {\frac{1}{\delta}}\right \rceil}   k_{f^{\alpha}} \delta^{s-1} .
\end{aligned}
\end{equation}
Consequently, we deduce
 $$\overline{\dim}_B\big(Graph(f^{\alpha})\big) =\varlimsup_{\delta \rightarrow 0} \frac{\log N_{\delta}(Graph(f^{\alpha}))}{- \log \delta}\le 2- s.$$
It is sufficient to prove the following bound for lower box dimension: $$\underline{\dim}_B\big(Graph(f^{\alpha})\big) \ge 2 - s.$$
We recall the self-referential equation $$ f^{\alpha}(x)= f(x) + \alpha_j\big(L_j^{-1}(x)\big)  \big[f^{\alpha}\big(L_j^{-1}(x)\big) - b\big(L_j^{-1}(x)\big)\big],$$ for every $x \in I_j $ and $j \in J.$
For $ x ,y \in I_j $ such that $|x-y| \le \delta,$ we obtain
\begin{equation*}
     \begin{aligned}
      |f^{\alpha}(x)- f^{\alpha}(y)| =& \Big| f(x)-f(y) + \alpha_j\big(L_j^{-1}(x)\big) ~ f^{\alpha}\big(L_j^{-1}(x)\big) - \alpha_j\big(L_j^{-1}(y)\big) ~ f^{\alpha}\big(L_j^{-1}(y)\big)\\& - \alpha_j\big(L_j^{-1}(x)\big) ~ b\big(L_j^{-1}(x)\big) + \alpha_j\big(L_j^{-1}(y)\big) ~ b\big(L_j^{-1}(y)\big) \Big|\\
      \ge & | f(x)-f(y)| - \|\alpha\|_{\infty} ~ \Big|f^{\alpha}\big(L_j^{-1}(x)\big) -  f^{\alpha}\big(L_j^{-1}(y)\big) \Big|\\& - \|\alpha\|_{\infty} ~ \Big|b\big(L_j^{-1}(x)\big) -   b\big(L_j^{-1}(y)\big) \Big|\\&- \big(\|b\|_{\infty}+\|f^{\alpha}\|_{\infty}\big) ~\Big|\alpha_j\big(L_j^{-1}(x)\big)-\alpha_j\big(L_j^{-1}(y)\big)
      \Big|
 \end{aligned}
 \end{equation*}
 With the help of Equation (\ref{Hypo}), we obtain
 \begin{equation*}
      \begin{aligned}
       |f^{\alpha}(x)- f^{\alpha}(y)|
       \ge & ~K_f | x-y|^{s}-\|\alpha\|_{\infty} ~ k_{f^\alpha}\Big|L_j^{-1}(x) -  L_j^{-1}(y) \Big|^{s}\\& - \|\alpha\|_{\infty}~ k_b \Big|L_j^{-1}(x) -   L_j^{-1}(y) \Big|^{s}\\&- \big(\|b\|_{\infty}+M\big) k_{\alpha} ~\Big|L_j^{-1}(x) -   L_j^{-1}(y) \Big|^{s}\\
       \ge & ~K_f | x-y|^{s}-\|\alpha\|_{\infty} ~ k_{f^\alpha} a^{-s} | x-y|^{s}\\& - \|\alpha\|_{\infty} ~ k_b a^{-s} | x-y|^{s}\\&- \big(\|b\|_{\infty}+M\big)a^{-s} k_{\alpha} ~|x-y|^{s}\\
      = & ~\Big(K_f-(k_{f^\alpha}+ k_b)\|\alpha\|_{\infty} a^{-s}-\big(\|b\|_{\infty}+M\big)a^{-s} k_{\alpha}\Big) | x-y|^{s}.
  \end{aligned}
  \end{equation*}
Let $K:=K_f-(k_{f^\alpha}+ k_b)\|\alpha\|_{\infty} a^{-s}-\big(\|b\|_{\infty}+M\big)a^{-s} k_{\alpha}.$ For $\delta= a^{m},$ we estimate 
       \begin{equation*}
            \begin{aligned}
             N_{\delta}(Graph(f^{\alpha})) &\ge  \sum_{i=1}^{a^{-m}} \max\Big\{1,\left \lceil {  a^{-m} R_{f^{\alpha}}[(i-1)\delta,i\delta] }\right \rceil \Big\}\\
             & \ge  \sum_{i=1}^{a^{-m} } \left \lceil {  a^{-m} R_{f^{\alpha}}[(i-1)\delta,i\delta] }\right \rceil \\ & \ge  \sum_{i=1}^{a^{-m} } \left \lceil {K a^{-m}   a^{m s}}\right \rceil\\
            & \ge   a^{-m}  a^{-m} K  a^{m s} \\
         & =   K a^{m( s-2)}   .
            \end{aligned}
            \end{equation*}
   Using the above bound for $N_{\delta}(Graph(f^{\alpha}))$, we obtain
\begin{equation*}
          \begin{aligned}
             \varliminf_{\delta \rightarrow 0}\frac{ \log\Big( N_{\delta}(Graph(f^{\alpha}))\Big)}{- \log (\delta)}
            & \ge \varliminf_{m \rightarrow \infty}\frac{ \log\Big(K  a^{m(s- 2)}   \Big) }{-m \log a}\\ & =
             2- s,
          \end{aligned}
          \end{equation*}
establishing the result.
\end{proof}

\begin{corollary}
If we consider the Bernstein polynomial as the base function, then for a Lipschitz $f$, we obtain a sequence of Bernstein $\alpha$-fractal functions (see \cite{CJN,Vijender} for details). For each $n\in \mathbb{N}$, let $G$ be the graph of the Bernstein $\alpha$-fractal function. Then under the hypothesis of Theorem \ref{mainthm}, we obtain $\dim_B(G)\leq 1$.   
\end{corollary}
\begin{remark}
In \cite{AGN2}, Nasim et al. computed the box dimension of $\alpha$-fractal function under certain condition. But for the H\"older exponent $s\in (0,1)$ the author has calculated the obvious upper bound as $2-s$. However, in this article, in Theorem \ref{mainthm}, we have calculated the exact estimation of the box dimension of $\alpha$-fractal function under suitable condition.
\end{remark}
\begin{theorem}\label{BBVL4}
    Let $f , \alpha_j (j \in J)$ and $b$ be  H\"older continuous with exponent $s$ such that $b(x_1)=f(x_1)$ and $b(x_N)=f(x_N)$. If 
    $  \|\alpha\|_{\mathcal{H}}< a^s $ with $a=\min\{a_j:j \in J\}$ then $$ 1 \le  \dim_H (Graph(f^{\alpha}) ) \le  2- s.$$
 \end{theorem}

\begin{proof}
We will proceed by defining a map $\Phi : Graph(f^{\alpha}) \rightarrow I $ by $\Phi((x,f(x)))=x.$ Then
$$|\Phi((x,f(x)))-\Phi((y,f(y)))| = |x-y| \le \|(x,f(x))-(y,f(y))\|_2 .$$ That is, $\Phi$ is a Lipschitz map. Using a properties of Hausdorff dimension (see \cite{Fal}), we have $\dim_H (\Phi(Graph(f^{\alpha}) )) \le \dim_H (Graph(f^{\alpha}) ).$ It is easy to check that the map $\Phi$ is onto. Hence we have $  \dim_H (Graph(f^{\alpha}) ) \ge \dim_H (I)=1 .$ We recall a well-known result, see \cite{Fal}, which relates the Hausdorff dimension and box dimension in the following sense: $$ \dim_H (C) \le \underline{\dim}_B (C) \le  \overline{\dim}_B (C) $$ for any bounded set $C \subset \mathbb{R}^n.$ Theorem \ref{BBVL3} and the first part of Theorem \ref{mainthm} yield the required upper bound for the Hausdorff dimension of the graph of fractal function $f^{\alpha}.$
\end{proof}

\begin{definition}
Let $ f:I \rightarrow \mathbb{R}$ be a function. For each partition $ P: t_0<t_1<t_2 < \dots <t_n $ of the interval $I, $ we define $$V(f,I)= \sup_P \sum_{i=1}^{n} |f(t_i)-f(t_{i-1})|,$$ where the supremum is taken over all partitions $P$ of the interval $I.$\\ If $V(f,I) < \infty,$ we say that $f$ is of bounded variation. The set of all functions of  bounded variation on $I$ will be denoted by $\mathcal{BV}(I)$.
We define a norm on $\mathcal{BV}(I)$ by $\|f\|_{\mathcal{BV}}:= |f(t_0)|+ V(f,I).$
Moreover, the space $\mathcal{BV}(I)$ is a Banach space with respect to this norm.
\end{definition}
Liang \cite{Liang1} proved that
\begin{theorem} \label{ET4}
If $f \in \mathcal{C}(I)\cap \mathcal{BV}(I)$, then $\dim_H(Graph(f))=\dim_B(Graph(f))=1.$
\end{theorem}

 The next remark is straightforward but useful for the upcoming theorem.
 \begin{remark}\label{Bound}
 Let $f$ be real-valued function on $I=[0,1].$ For $c,d \in \mathbb{R},$ we define a function $g(x)=f(cx+d) $ on a suitable domain. If $f$ is of bounded variation on $I$ then $g$ is also of bounded variation on its domain.
 
 \end{remark}
 We present the following remark for the sake of independent interest.
 \begin{remark}
 We know that $f^{\alpha}$ satisfies the self-referential equation
 $$ f^{\alpha}(x)= f(x)+\alpha_j(L_j^{-1}(x)).(f^{\alpha}- b)\big(L_j^{-1}(x)\big)~~~~\forall~~ x \in I_j,~~ j \in J.$$
 The self-referential equation may also be written in the following manner
 $$ \alpha_j(L_j^{-1}(x))~b\big(L_j^{-1}(x)\big)= f(x)- f^{\alpha}(x)+\alpha_j(L_j^{-1}(x))~f^{\alpha}\big(L_j^{-1}(x)\big) ~~~~\forall~~ x \in I_j,~~ j \in J.$$
 Further, we assume $f,f^{\alpha}$ and $\alpha_j$ ($j \in J$) be of bounded variation  with $\alpha_j> 0 ~\text{or} < 0$ on $I$. Using algebra of bounded variation functions (see \cite{Gordon}), one concludes that $b$ is of bounded variation. 
 \end{remark}
 
 The following theorem is a generalization of \cite[Theorem $4.8$]{VV}. However, we present the proof for reader's convenience.
 \begin{theorem} \label{bounded2}
 Let $f \in \mathcal{BV}(I).$ Suppose that $\triangle=\{x_1,x_2\dots,x_N:x_1 <x_2< \dots <x_N \}$ is a partition of $I ,$ $b \in \mathcal{BV}(I)$ satisfying $b(x_1)=f(x_1),~b(x_N)=f(x_N)$, and $\alpha_j ~(j \in J)$ are functions in $\mathcal{BV}(I)$ with $\|\alpha\|_{\mathcal{BV}}< \frac{1}{2(N-1)}.$ Then, the fractal perturbation $ f^{\alpha}$ corresponding to $f$ is of bounded variation on $I$.

 \end{theorem}
 \begin{proof}
 Let $\mathcal{BV}_*(I)= \{g \in \mathcal{BV}(I):g(x_1)=f(x_1), ~g(x_N)=f(x_N)\}.$ We may see (using the definition of closed set) that $\mathcal{BV}_*(I)$ is a closed subset of $\mathcal{BV}(I)$. Being a close subset of Banach space $\mathcal{BV}(I)$, the space $\mathcal{BV}_*(I)$ is a complete metric space when endowed with metric induced by norm $\|g\|_{\mathcal{BV}}:= |g(x_1)|+ V(g,I).$
 Define the RB operator $T: \mathcal{BV}_*(I) \rightarrow \mathcal{BV}_*(I) $ by $$ (Tg)(x)= f(x) + \alpha_j\big(L_j^{-1}(x)\big)  \big[g\big(L_j^{-1}(x)\big) - b\big(L_j^{-1}(x)\big)\big],$$ for every $x \in I_j$ and $j \in J.$ As done in previous theorems we note that  $T$ is well-defined. Let $ P: t_0<t_1<t_2 < \dots <t_m $ be a partition of the interval $I_j,$ where $m \in \mathbb{N}.$ Consider 
 \begin{equation*}
 \begin{aligned}
   \Big|(Tg-Th)(t_i)- (Tg-Th)(t_{i-1})\Big|=& ~ \Big| \alpha_j\big(L_j^{-1}(t_i)\big)(g -h)\big(L_j^{-1}(t_i)\big)\\ &  -  \alpha_j\big(L_j^{-1}(t_{i-1})\big)(g -h)\big(L_j^{-1}(t_{i-1})\big)\Big|\\
    \le & ~ | \alpha_j\big(L_j^{-1}(t_i)\big)| \Big|(g -h)\big(L_j^{-1}(t_i)\big) \\&-  (g -h)\big(L_j^{-1}(t_{i-1})\big)\Big|+\Big|(g -h)\big(L_j^{-1}(t_{i-1})\big)\Big|\\
        &. \Big|\alpha_j\big(L_j^{-1}(t_i)\big) -\alpha_j\big(L_j^{-1}(t_{i-1})\big) \Big|\\
   \le & ~ \|\alpha\|_{\infty} \Big|(g -h)\big(L_j^{-1}(t_i)\big) -  (g -h)\big(L_j^{-1}(t_{i-1})\big)\Big| \\
   & +   \|g-h\|_{\infty} \Big|\alpha_j\big(L_j^{-1}(t_i)\big) -\alpha_j\big(L_j^{-1}(t_{i-1})\big) \Big|.
 \end{aligned} 
 \end{equation*}
 Summing over $i=1 $ to $m,$ we have
 \begin{equation*}
 \begin{aligned}
       \sum_{i=1}^{m}&\big|(Tg-Th)(t_i)- (Tg-Th)(t_{i-1})\big| \\ ~~\le & ~\|\alpha\|_{\infty} \sum_{i=1}^{m} \big|(g -h)\big(L_j^{-1}(t_i)\big) -  (g -h)\big(L_j^{-1}(t_{i-1})\big)\big| \\
      &+ \|g-h\|_{\infty} \sum_{i=1}^{m} \Big|\alpha_j\big(L_j^{-1}(t_i)\big) -\alpha_j\big(L_j^{-1}(t_{i-1})\big) \Big|\\
      \le &~ \|\alpha\|_{\infty} \|g-h\|_{\mathcal{BV}}+ \|g-h\|_{\infty} \|\alpha\|_{\mathcal{BV}}\\
      \le & ~ \|\alpha\|_{\mathcal{BV}} \Big( \|g-h\|_{\mathcal{BV}}+\|g-h\|_{\infty}\Big)\\
      \le & ~2 \|\alpha\|_{\mathcal{BV}}~ \|g-h\|_{\mathcal{BV}}.
 \end{aligned} 
 \end{equation*}
 The above inequality holds for any partition of $I_j$. Therefore, one gets $$ \|Tg-Th\|_{\mathcal{BV}} \le 2 (N-1) \|\alpha\|_{\mathcal{BV}} \|g-h\|_{\mathcal{BV}}.$$
 Since $\|\alpha\|_{\mathcal{BV}} < \frac{1}{2(N-1)},$ $T$ is a contraction on the complete metric space $\mathcal{BV}_*(I).$ Applying the Banach fixed point theorem we have a unique fixed point $f^{\alpha}$ of $T.$ Moreover, the fixed point $f^{\alpha}$ of $T$ satisfies the self-referential equation, that is, $$ f^{\alpha}(x)= f(x) + \alpha_j\big(L_j^{-1}(x)\big)  \big[f^{\alpha}\big(L_j^{-1}(x)\big) - b\big(L_j^{-1}(x)\big)\big],$$ for every $ x \in I_j $ and $j \in J.$

 \end{proof}

 \begin{theorem}
 Let the germ function $f$ and the parameter $b$ be continuous functions of bounded variation. Suppose $\alpha_j$ ($j \in J$) are functions of bounded variation with $\|\alpha\|_{\mathcal{BV}}< \frac{1}{2(N-1)}.$ Then $\dim_H(Graph(f^{\alpha}))=\dim_B(Graph(f^{\alpha}))=1.$
 \end{theorem}
 \begin{proof}
 Theorem \ref{ET4} and Theorem \ref{bounded2} produce the result. 
 \end{proof}
 We shall denote by $\mathcal{AC}(I)$ the Banach space of all absolutely continuous functions on $I$ with its usual norm (denoted by $\|.\|_{\mathcal{AC}}$). 
 
 \begin{theorem} \label{bounded3}
  Let $f \in \mathcal{AC}(I) .$ Suppose that $\triangle=\{ x_1,x_2,\dots,x_N:x_1 <x_2< \dots <x_N \}$ is a partition of $I ,$ $b \in \mathcal{AC}(I)$ satisfying $b(x_1)=f(x_1),~b(x_N)=f(x_N)$ , and $\alpha_j ~(j \in J)$ are functions in $\mathcal{AC}(I)$ with $\|\alpha\|_{\mathcal{AC}} < \frac{a}{2(N-1)},$ where $a=\min\{a_j: j \in J\}.$ Then, the fractal perturbation $ f^{\alpha}$ corresponding to $f$ is absolutely continuous on $I$.

  \end{theorem}
  \begin{proof}
  Let $\mathcal{AC}_*(I)= \{g \in \mathcal{AC}(I):g(x_1)=f(x_1), ~g(x_N)=f(x_N)\}.$ We observe (using the sequential definition of a closed set) that $\mathcal{AC}_*(I)$ is a closed subset of $\mathcal{AC}(I)$. Since $\mathcal{AC}(I)$ endowed with $\|g\|_{\mathcal{AC}}:= |g(x_1)|+ \int_{x_1}^{x_N} |g'(x)|dx$ is a complete normed linear space, the set $\mathcal{AC}_*(I)$ is a complete metric space when equipped with metric induced by aforesaid norm.
  Define the RB operator $T: \mathcal{AC}_*(I) \rightarrow \mathcal{AC}_*(I) $ by $$ (Tg)(x)= f(x) + \alpha_j\big(L_j^{-1}(x)\big)  \big[g\big(L_j^{-1}(x)\big) - b\big(L_j^{-1}(x)\big)\big],$$ for every $x \in I_j$ and $j \in J.$ We note that the conditions on $f$ and $b$ dictate the function $T$ to be well-defined. Consider 
  \begin{equation*}
  \begin{aligned}
    \int_{L_j(x_1)}^{L_j(x_N)}|(Tg-Th)'(x)|dx \le & ~\frac{1}{a_j} \int_{L_j(x_1)}^{L_j(x_N)} | \alpha_j'\big(L_j^{-1}(x)\big)(g -h)\big(L_j^{-1}(x)\big) |dx \\
    &+ \frac{1}{a_j} \int_{L_j(x_1)}^{L_j(x_N)} |\alpha_j\big(L_j^{-1}(x)\big)(g -h)'\big(L_j^{-1}(x)\big)| dx \\
     = & ~\frac{1}{a_j} \int_{x_1}^{x_N} | \alpha_j'(y)(g -h)(y)|dy \\
         &+ \frac{1}{a_j} \int_{x_1}^{x_N} |\alpha_j(y)(g -h)'(y)| dy\\
    \le & ~\frac{\|g-h\|_{\infty}}{a_j} \int_{x_1}^{x_N} | \alpha_j'(y)|dy \\
             &+ \frac{\|\alpha_j\|_{\infty}}{a_j} \int_{x_1}^{x_N} |(g -h)'(y)| dy.
  \end{aligned} 
  \end{equation*}
  Summing over $j=1 $ to $N-1,$ we have
  \begin{equation*}
  \begin{aligned}
       \sum_{j=1}^{N-1}\int_{L_j(x_1)}^{L_j(x_N)}|(Tg-Th)'(x)|dx \le &  ~ \frac{2(N-1)\|\alpha\|_{\mathcal{AC}}}{a}  \|g-h\|_{\mathcal{AC}}.
  \end{aligned} 
  \end{equation*}
  Therefore, one gets $$ \|Tg-Th\|_{\mathcal{AC}} \le \frac{2(N-1)\|\alpha\|_{\mathcal{AC}}}{a}  \|g-h\|_{\mathcal{AC}}.$$
  Since $\|\alpha\|_{\mathcal{AC}} < \frac{a}{2(N-1)},$ we deduce that $T$ is a contraction on the complete metric space $\mathcal{AC}_*(I).$ Moreover, the fixed point $f^{\alpha}$ of $T$ satisfies the self-referential equation, that is, $$ f^{\alpha}(x)= f(x) + \alpha_j\big(L_j^{-1}(x)\big)  \big[f^{\alpha}\big(L_j^{-1}(x)\big) - b\big(L_j^{-1}(x)\big)\big],$$ for every $ x \in I_j $ and $j \in J.$

  \end{proof}
  
  Combining Theorem \ref{ET4} and Theorem \ref{bounded2}, one can immediately deduce the following.
  \begin{theorem}
  Let the germ function $f$ and the parameter $b$ be absolutely continuous functions. Suppose $\alpha_j$ ($j \in J$) are absolutely continuous functions with $\|\alpha\|_{\mathcal{AC}} < \frac{a}{2(N-1)}.$ Then  $\dim_H(Graph(f^{\alpha}))=\dim_B(Graph(f^{\alpha}))=1.$
  \end{theorem}

\subsection*{Acknowledgements}
  The first author thanks Dr. P. Viswanathan for his suggestions, support and encouragement during preparation of the manuscript.

\bibliographystyle{amsplain}
\bibliography{Dimension.bib}

\providecommand{\bysame}{\leavevmode\hbox to3em{\hrulefill}\thinspace}
\providecommand{\MR}{\relax\ifhmode\unskip\space\fi MR }
\providecommand{\MRhref}[2]{%
  \href{http://www.ams.org/mathscinet-getitem?mr=#1}{#2}
}
\providecommand{\href}[2]{#2}
\begin{thebibliography}{10}

\bibitem{AGN2}
Md.~N. Akhtar, M.~G.~P. Prasad, and M.~A. Navascu\'es, \emph{Box dimensions of
  {$\alpha$}-fractal functions}, Fractals \textbf{24} (2016), no.~3, 1650037,
  13.

\bibitem{Bandt}
C.~Bandt and S.~Graf, \emph{Self-similar sets. {VII}. {A} characterization of
  self-similar fractals with positive {H}ausdorff measure}, Proc. Amer. Math.
  Soc. \textbf{114} (1992), no.~4, 995--1001.

\bibitem{MF1}
M.~F. Barnsley, \emph{Fractal functions and interpolation}, Constr. Approx.
  \textbf{2} (1986), no.~4, 303--329.

\bibitem{MF2}
\bysame, \emph{Fractals {E}verywhere}, Academic Press, Inc., Boston, MA, 1988.

\bibitem{MF6}
M.~F. Barnsley, J.~Elton, D.~Hardin, and P.~Massopust, \emph{Hidden variable
  fractal interpolation functions}, SIAM J. Math. Anal. \textbf{20} (1989),
  no.~5, 1218--1242.

\bibitem{MF3}
M.~F. Barnsley and P.~R. Massopust, \emph{Bilinear fractal interpolation and
  box dimension}, J. Approx. Theory \textbf{192} (2015), 362--378.

\bibitem{ACar}
A.~Carvalho, \emph{Box dimension, oscillation and smoothness in function
  spaces}, J. Funct. Spaces Appl. \textbf{3} (2005), no.~3, 287--320.

\bibitem{CJN}
A.~K.~B. Chand, S.~Jha, and M.~A. Navascu\'es, \emph{Kantorovich-{B}ernstein
  $\alpha$-fractal functions in $\mathcal{L}^p$ spaces}, Quaest. Math.
  \textbf{43} (2020), no.~2, 227--241.

\bibitem{ADBJ}
A.~Deliu and B.~Jawerth, \emph{Geometrical dimension versus smoothness},
  Constr. Approx. \textbf{8} (1992), no.~2, 211--222.

\bibitem{Fal}
K.~Falconer, \emph{Fractal {G}eometry}, second ed., John Wiley \& Sons, Inc.,
  Hoboken, NJ, 2003, Mathematical {F}oundations and {A}pplications.

\bibitem{Fal2}
K.~J. Falconer, \emph{Dimensions and measures of quasi self-similar sets},
  Proc. Amer. Math. Soc. \textbf{106} (1989), no.~2, 543--554.

\bibitem{Lau1}
A.~H. Fan and K.~S. Lau, \emph{Iterated function system and {R}uelle operator},
  J. Math. Anal. Appl. \textbf{231} (1999), no.~2, 319--344.

\bibitem{Gordon}
R.~A. Gordon, \emph{Real {A}nalysis: A first course}, second ed., Addison
  Wesley, 2001.

\bibitem{Hardin}
D.~P. Hardin and P.~R. Massopust, \emph{The capacity for a class of fractal
  functions}, Comm. Math. Phys. \textbf{105} (1986), no.~3, 455--460.

\bibitem{HM}
\bysame, \emph{Fractal interpolation functions from {$\bold R^n$} into {$\bold
  R^m$} and their projections}, Z. Anal. Anwendungen \textbf{12} (1993), no.~3,
  535--548.

\bibitem{Lau2}
K.~S. Lau, H.~Rao, and Y.~L. Ye, \emph{Corrigendum: ``{I}terated function
  system and {R}uelle operator'' [{J}. {M}ath. {A}nal. {A}ppl. {\bf 231}
  (1999), no. 2, 319--344; {MR}1669203 (2001a:37013)] by {L}au and {A}. {H}.
  {F}an}, J. Math. Anal. Appl. \textbf{262} (2001), no.~1, 446--451.

\bibitem{Liang1}
Y.~S. Liang, \emph{Box dimensions of {R}iemann-{L}iouville fractional integrals
  of continuous functions of bounded variation}, Nonlinear Anal. \textbf{72}
  (2010), no.~11, 4304--4306.

\bibitem{M2}
P.~R. Massopust, \emph{Fractal {F}unctions, {F}ractal {S}urfaces, and
  {W}avelets}, Academic Press, Inc., San Diego, CA, 1994.

\bibitem{M1}
\bysame, \emph{Interpolation and {A}pproximation with {S}plines and
  {F}ractals}, Oxford University Press, Oxford, 2010.

\bibitem{N2}
M.~A. Navascu\'es, \emph{Fractal polynomial interpolation}, Z. Anal.
  Anwendungen \textbf{24} (2005), no.~2, 401--418.

\bibitem{N1}
\bysame, \emph{Fractal approximation}, Complex Anal. Oper. Theory \textbf{4}
  (2010), no.~4, 953--974.

\bibitem{Peres}
Y.~Peres, M.~Rams, K.~Simon, and B.~Solomyak, \emph{Equivalence of positive
  {H}ausdorff measure and the open set condition for self-conformal sets},
  Proc. Amer. Math. Soc. \textbf{129} (2001), no.~9, 2689--2699.

\bibitem{RSY}
H.~J. Ruan, W.~Y. Su, and K.~Yao, \emph{Box dimension and fractional integral
  of linear fractal interpolation functions}, J. Approx. Theory \textbf{161}
  (2009), no.~1, 187--197.

\bibitem{Schief}
A.~Schief, \emph{Separation properties for self-similar sets}, Proc. Amer.
  Math. Soc. \textbf{122} (1994), no.~1, 111--115.

\bibitem{VV}
S.~Verma and P.~Viswanathan, \emph{A revisit to {$\alpha$}-fractal function and
  box dimension of its graph}, Fractals \textbf{27} (2019), no.~6, 1950090, 15.

\bibitem{VV1}
\bysame, \emph{A fractal operator associated with bivariate fractal
  interpolation functions on rectangular grids}, Results Math. \textbf{75}
  (2020), no.~1, Paper No. 28, 26.

\bibitem{Vijender2}
N.~Vijender, \emph{Approximation by hidden variable fractal functions: a
  sequential approach}, Results Math. \textbf{74} (2019), no.~4, Paper No. 192,
  23.

\bibitem{Vijender}
\bysame, \emph{Bernstein fractal trigonometric approximation}, Acta Appl. Math.
  \textbf{159} (2019), 11--27.

\bibitem{WY}
H.~Y. Wang and J.~S. Yu, \emph{Fractal interpolation functions with variable
  parameters and their analytical properties}, J. Approx. Theory \textbf{175}
  (2013), 1--18.

\bibitem{Lau3}
Yuan-Ling Ye, \emph{Separation properties for self-conformal sets}, Studia
  Math. \textbf{152} (2002), no.~1, 33--44.

\end{thebibliography}

\end{document}